\documentclass[12pt]{amsart}

\usepackage{amssymb}
\usepackage{graphicx}
\usepackage{ifthen}
\usepackage{mathrsfs}
\usepackage{pict2e}
\usepackage{xargs}
\usepackage{xspace}

\newcommand{\definedterm}[1]{\emph{#1}}

\newcommand{\Bairespace}[1][]{
  \ifthenelse{\equal{#1}{}}{\functions{\N}{\N}}{\functions{#1}{\N}}
}

\newcommand{\Bairetree}[1][]{
  \ifthenelse{\equal{#1}{}}{\functions{<\N}{\N}}{\functions{#1}{\N}}
}

\newcommand{\calF}{\mathscr{F}}

\newcommand{\Cantorspace}[1][\N]{
  \functions{#1}{2}
}

\newcommand{\cardinality}[1]{|#1|}

\newcommand{\composition}{\circ}
\newcommandx{\concatenation}[2][1 = undefined, 2 = undefined]{
  \ifthenelse{\equal{#1}{undefined}}{{}\smallfrown}{
    \ifthenelse{\equal{#2}{undefined}}{\bigoplus #1}{\bigoplus_{#1} #2}
  }
}

\newcommand{\cost}[2]{C_{#1}(#2)}

\newcommand{\domain}[1]{\mathrm{dom}(#1)}

\newcommand{\equivalenceclass}[2]{[#1]_{#2}}

\newcommand{\evaluation}{\text{eval}}

\newcommand{\F}[1]{\mathbb{F}_{#1}}

\newcommand{\from}{\colon}
\newcommandx{\functions}[2]{#2^{#1}}

\newcommand{\graph}[1]{\mathrm{graph}(#1)}

\newcommand{\horizontalsection}[2]{#1^{#2}}

\newcommand{\id}{\mathrm{id}}
\newcommand{\image}[2]{#1(#2)}
\newcommandx{\intersection}[2][1 = undefined, 2 = undefined]{
  \ifthenelse{\equal{#1}{undefined}}{\cap}{
    \ifthenelse{\equal{#2}{undefined}}{\bigcap #1}{\bigcap_{#1} #2}
  }
}
\newcommand{\inverse}[1]{#1^{-1}}
\newcommandx{\join}[2][1 = undefined, 2 = undefined]{
  \ifthenelse{\equal{#1}{undefined}}{\vee}{
    \ifthenelse{\equal{#2}{undefined}}{\bigvee #1}{\bigvee_{#1} #2}
  }
}

\newcommand{\mathand}{\text{ and }}

\newcommand{\N}{\mathbb{N}}

\newcommand{\orbitequivalencerelation}[2]{E_{#1}^{#2}}
\newcommand{\pair}[2]{(#1, #2)}
\newcommand{\preimage}[2]{#1^{-1}(#2)}

\newcommand{\projection}[1]{\mathrm{proj}_{#1}}

\newcommand{\R}{\mathbb{R}}

\renewcommand{\restriction}[2]{#1 \upharpoonright #2}

\newcommand{\saturation}[2]{[#1]_{#2}}
\newcommandx{\sequence}[2][2 = undefined]{
  \ifthenelse{\equal{#2}{undefined}}{(#1)}{
    (#1)_{#2}
  }
}
\newcommandx{\set}[2][2 = undefined]{
  \ifthenelse{\equal{#2}{undefined}}{\{ #1 \}}{
    \{ #1 \suchthat #2 \}
  }
}
\newcommandx{\sets}[4][3 = undefined, 4 = undefined]{
  \ifthenelse{\equal{#3}{undefined}}{[#2]^{#1}}{
    \ifthenelse{\equal{#4}{undefined}}{[#2]^{#1}_{#3}}{[#2]^{#1}_{#3 / #4}}
  }
}

\newcommand{\sigmaclass}[1]{\sigma(#1)}
\newcommandx{\Sigmaclass}[2][1=,2=]{
  \ifthenelse{\equal{#2}{}}{\mathbf{\Sigma}_{#1}}{\mathbf{\Sigma}^{#1}_{#2}}
}
\newcommand{\singleton}[1]{\set{#1}}

\newcommand{\SL}[2]{\mathrm{SL}_{#1}(#2)}
\renewcommand{\square}[1]{I(#1)}
\newcommand{\suchthat}{\mid}
\newcommand{\support}[1]{\mathrm{supp}(#1)}
\newcommand{\T}[1][]{\mathbb{T}^{#1}}

\newcommand{\uniformmetric}[1]{d_{#1}}
\newcommandx{\union}[2][1 = undefined, 2 = undefined]{
  \ifthenelse{\equal{#1}{undefined}}{\cup}{
    \ifthenelse{\equal{#2}{undefined}}{\bigcup #1}{\bigcup_{#1} #2}
  }
}
\newcommand{\verticalsection}[2]{#1_{#2}}
\newcommand{\Z}{\mathbb{Z}}

\newcommand{\Baire}{Baire\xspace}
\newcommand{\Borel}{Bor\-el\xspace}

\newcommand{\Dye}{Dye\xspace}
\newcommand{\Farrell}{Far\-rell\xspace}
\newcommand{\Feldman}{Feld\-man\xspace}
\newcommand{\Gaboriau}{Gab\-or\-i\-au\xspace}

\newcommand{\Hopf}{Hopf\xspace}

\newcommand{\Jankov}{Jan\-kov\xspace}
\newcommand{\Kechris}{Kech\-ris\xspace}

\newcommand{\Kuratowski}{Kur\-at\-ow\-ski\xspace}

\newcommand{\Lusin}{Lu\-sin\xspace}
\newcommand{\Miller}{Mil\-ler\xspace}
\newcommand{\Moore}{Moore\xspace}
\newcommand{\Mycielski}{My\-ciel\-ski\xspace}
\newcommand{\Novikov}{No\-vik\-ov\xspace}
\newcommand{\Ornstein}{Orn\-stein\xspace}
\newcommand{\Pichot}{Pi\-chot\xspace}
\newcommand{\Polish}{Po\-lish\xspace}

\newcommand{\Sierpinski}{Sier\-pi\'{n}\-ski\xspace}

\newcommand{\Ulam}{U\-lam\xspace}
\newcommand{\Varadarajan}{Var\-ad\-ar\-aj\-an\xspace}
\newcommand{\vonNeumann}{von Neu\-mann\xspace}
\newcommand{\Weiss}{Weiss\xspace}

\newenvironment{lemmaproof}{
  
  \begin{proof}
}{\end{proof}}

\newenvironment{propositionproof}{
  
  \begin{proof}
}{\end{proof}}

\newenvironment{theoremproof}{
  
  \begin{proof}
}{\end{proof}}

\newtheorem{lemma}{Lemma}[section]
\newtheorem{proposition}[lemma]{Proposition}
\newtheorem{theorem}[lemma]{Theorem}

\newtheorem*{introtheorem}{Theorem}

\theoremstyle{definition}

\newtheorem{remark}[lemma]{Remark}
\newtheorem*{acknowledgments}{Acknowledgments}

\begin{document}

\begin{abstract}
  Suppose that $X$ is a \Polish space, $E$ is a countable \Borel equivalence
  relation on $X$, and $\mu$ is an $E$-invariant \Borel probability measure on $X$. We
  consider the circumstances under which for every countable non-abelian free group
  $\Gamma$, there is a \Borel sequence $\sequence{\cdot_r}[r \in \R]$ of free actions of
  $\Gamma$ on $X$, generating subequivalence relations $E_r$ of $E$ with respect to
  which $\mu$ is ergodic, with the further property that $\sequence{E_r}[r \in \R]$ is an
  increasing sequence of relations which are pairwise incomparable under
  $\mu$-reducibility. In particular, we show that if $E$ satisfies a natural separability
  condition, then this is the case as long as there exists a free \Borel action of a countable
  non-abelian free group on $X$, generating a subequivalence relation of $E$ with respect
  to which $\mu$ is ergodic.
\end{abstract}

\author[C.T. Conley]{Clinton T. Conley}

\address{
  Clinton T. Conley \\
  Department of Mathematical Sciences \\
  Carnegie Mellon University \\
  Pittsburgh, PA 15213 \\
  USA
}

\email{clintonc@andrew.cmu.edu}

\urladdr{
  http://www.math.cmu.edu/math/faculty/Conley
}

\author[B.D. Miller]{Benjamin D. Miller}

\address{
  Benjamin D. Miller \\
  Kurt G\"{o}del Research Center for Mathematical Logic \\
  Universit\"{a}t Wien \\
  W\"{a}hringer Stra{\ss}e 25 \\
  1090 Wien \\
  Austria
 }

\email{benjamin.miller@univie.ac.at}

\urladdr{
  http://www.logic.univie.ac.at/benjamin.miller
}

\thanks{The authors were supported in part by FWF Grant P28153 and SFB Grant 878.}

  

\title{Incomparable actions of free groups}

\maketitle

\section*{Introduction}

Suppose that $E$ and $F$ are equivalence relations on $X$ and $Y$. A \definedterm
{homomorphism} from $E$ to $F$ is a function $\phi \from X \to Y$ sending $E$-equivalent points
to $F$-equivalent points, a \definedterm{reduction} of $E$ to $F$ is a homomorphism sending
$E$-inequivalent points to $F$-inequivalent points, an \definedterm{embedding} of $E$ into
$F$ is an injective reduction, and an \definedterm{isomorphism} of $E$ with $F$ is a bijective 
reduction. \Borel reducibility plays a central role in the descriptive set-theoretic study of equivalence
relations, whereas measurable isomorphism figures more prominently in the ergodic-theoretic
study of equivalence relations.

We refer the reader to \cite{Kechris} for basic descriptive-set-theoretic definitions
and notation.

The \definedterm{orbit equivalence relation} generated by an action of a group $\Gamma$ on a
set $X$ is given by $x \mathrel{\orbitequivalencerelation{\Gamma}{X}} y
\iff \exists \gamma \in \Gamma \ \gamma \cdot x = y$. Actions of groups $\Gamma$ and
$\Delta$ on \Polish spaces $X$ and $Y$ equipped with \Borel measures $\mu$ and
$\nu$ are \definedterm{orbit equivalent} if there is a measure-preserving \Borel isomorphism
of the restrictions of their orbit equivalence relations to conull \Borel sets.

Following the usual abuse of language, we say that an equivalence relation is \definedterm
{countable} if its classes are all countable, and \definedterm{finite} if its classes are all finite.
Suppose now that $E$ is a countable \Borel equivalence relation on $X$. We say that a \Borel 
probability measure $\mu$ on $X$ is \definedterm{$E$-ergodic} if every $E$-invariant \Borel
set is $\mu$-null or $\mu$-conull. We say that $\mu$ is \definedterm{$E$-invariant} if $\mu
(B) = \mu(\image{T}{B})$, for every \Borel set $B \subseteq X$ and every \Borel automorphism
$T \from X \to X$ whose graph is contained in $E$. And we say that $\mu$ is \definedterm
{$E$-quasi-invariant} if the family of $\mu$-null \Borel subsets of $X$ is closed under
$E$-saturation. As a result of \Feldman-\Moore (see \cite[Theorem 1]{FeldmanMoore})
ensures that $E$ is generated by a \Borel action of a countable discrete group, it follows that
every $E$-invariant \Borel probability measure is also $E$-quasi-invariant. 

Results of \Dye and \Ornstein-\Weiss ensure that all \Borel actions of countable amenable
groups equipped with suitably non-trivial ergodic invariant \Borel probability measures are
orbit equivalent (see, for example, \cite[Theorem 10.7]{KechrisMiller}). In contrast, the primary
result of \cite{GaboriauPopa} ensures that for every countable non-abelian free group
$\Gamma$, there are uncountably many \Borel-probability-measure-preserving ergodic free
\Borel actions of $\Gamma$ on a \Polish space which are pairwise non-orbit-equivalent. A
number of simplifications of the proof and strengthenings of the result have subsequently
appeared (see, for example, \cite{Tornquist, Hjorth}).

We say that $E$ is \definedterm{$\mu$-reducible} to $F$ if there is a $\mu$-conull \Borel set
on which there is a \Borel reduction of $E$ to $F$. Similarly, we say that $E$ is \definedterm
{$\mu$-somewhere reducible} to $F$ if there is a $\mu$-positive \Borel set on which there is
a \Borel reduction of $E$ to $F$. As the \Lusin-\Novikov uniformization theorem for \Borel
subsets of the plane with countable vertical sections (see, for example, \cite[Theorem
18.10]{Kechris}) ensures that every countable \Borel equivalence relation is \Borel reducible
to its restriction to any \Borel set intersecting every equivalence class, it follows that
when $\mu$ is $E$-ergodic, these two notions are equivalent.

In \cite{ConleyMiller}, a variety of results were established concerning the nature of
countable \Borel equivalence relations at the base of the measure reducibility hierarchy,
using arguments substantially simpler than those previously appearing. While the analogs
of many of these results for equivalence relations generated by free \Borel actions of
countable non-abelian free groups trivially follow, the analog corresponding to the main
results of \cite{GaboriauPopa, Hjorth} also requires a generalization of a stratification
result utilized in \cite{Hjorth, ConleyMiller}. Here we establish the latter, and incorporate it
into ideas from \cite{ConleyMiller} to obtain the former.

We say that a countable \Borel equivalence relation is \definedterm{hyperfinite} if it is the
union of an increasing sequence $\sequence{F_n}[n \in \N]$ of finite \Borel subequivalence
relations, \definedterm{$\mu$-hyperfinite} if it is hyperfinite on a $\mu$-conull \Borel set,
and \definedterm{$\mu$-nowhere hyperfinite} if there is no $\mu$-positive \Borel
set on which it is hyperfinite. We say that a countable \Borel equivalence relation $F$ on a
\Polish space $Y$ is \definedterm{projectively separable} if whenever $X$ is a
\Polish space, $E$ is a countable \Borel equivalence relation on $X$, and $\mu$
is a \Borel probability measure on $X$ for which $E$ is $\mu$-nowhere hyperfinite, the
pseudo-metric $\uniformmetric{\mu}(\phi, \psi) = \mu(\set{x \in X}[\phi(x) \neq \psi(x)])$ on the
space of all countable-to-one \Borel homomorphisms $\phi \from B \to Y$ from
$\restriction{E}{B}$ to $F$, where $B$ ranges over all \Borel subsets of $X$, is separable.

In \cite{ConleyMiller}, it is shown that the family of such equivalence relations includes the
orbit equivalence relation induced by the usual action of $\SL{2}{\Z}$ on $\T[2]$, and is
closed downward under both \Borel reducibility and \Borel subequivalence relations. 
Moreover, it is shown that such relations possess many of the exotic properties of 
countable \Borel equivalence relations which were previously known to hold only of
relations relatively high in the \Borel reducibility hierarchy. Our primary result here
strengthens one such theorem.

\begin{introtheorem} \label{introduction:mainresult}
  Suppose that $X$ is a \Polish space, $E$ is a projectively separable countable
  \Borel equivalence relation on $X$, $\mu$ is an $E$-invariant \Borel probability measure
  on $X$, and there is a free \Borel action of a countable non-abelian free group on $X$
  generating a subequivalence relation of $E$ with respect to which $\mu$ is ergodic. Then
  for every countable non-abelian free group $\Gamma$, there is a \Borel sequence
  $\sequence{\cdot_r}[r \in \R]$ of free actions of $\Gamma$ on $X$, generating 
  subequivalence relations $E_r$ of $E$ with respect to which $\mu$ is ergodic, with the
  further property that $\sequence{E_r}[r \in \R]$ is an increasing sequence of relations which
  are pairwise incomparable under $\mu$-reducibility.
\end{introtheorem}

The main result of \cite{GaboriauLyons} shows that every countable group has a free
\Borel action whose orbit equivalence relation satisfies the piece of our hypothesis
concerning the existence of an appropriate action of an appropriate free group. It is
an open question as to whether the failure of this requirement on a $\mu$-conull \Borel set,
or even its weakening in which ergodicity of the orbit equivalence relation is not
required, is equivalent to $\mu$-hyperfiniteness (see \cite[\S28]{KechrisMiller}).

Although the existence of incomparable orbit equivalence relations is primarily of interest in
the presence of \Borel probability measures which are both ergodic and invariant, we also
establish analogous results in which these assumptions are omitted (see Theorem 
\ref{antichains:actions}).

In \S\ref{ergodicgenerator}, we characterize the existence of subequivalence 
relations induced by free \Borel actions of countable non-abelian free groups in which one of
the generators acts ergodically. In \S\ref{stratification}, we use this to obtain new stratification 
results for free Borel actions of countable non-abelian free groups. And in \S\ref{incomparable},
we prove our primary results.

\section{Free actions with an ergodic generator} \label{ergodicgenerator}

Given a \Polish space $X$, a countable \Borel equivalence relation
$E$ on $X$, and a \Borel function $\phi \from X \to \N \union \set{\infty}$, one can
ask whether there is a \Borel function $\psi \from X \to X$, whose graph is
contained in $E$, such that $\cardinality{\preimage{\psi}{x}} = \phi(x)$ for
all $x \in X$. The uniformization theorem for \Borel subsets of the plane with
countable vertical sections ensures that if $\mu$ is an $E$-invariant \Borel probability
measure on $X$, then the existence of such a function necessitates that $\int \phi(x)
\ d\mu(x) = 1$. We begin this section by noting that if $\mu$ is also
$E$-ergodic, then this necessary condition is also sufficient (off of a $\mu$-null \Borel set).

The \definedterm{support} of $\phi \from X \to \N \union \set{\infty}$
is given by $\support{\phi} = X \setminus \preimage{\phi}{0}$. We say that
$\psi \from X \to Y \subseteq X$ is a \definedterm{retraction} if 
$\psi(y) = y$ for all $y \in Y$.

\begin{proposition} \label{ergodicgenerator:manytoone}
  Suppose that $X$ is a \Polish space, $E$ is a countable \Borel
  equivalence relation on $X$, $\phi \from X \to \N \union \set{\infty}$ is
  \Borel, and $\mu$ is an $E$-ergodic $E$-invariant \Borel probability
  measure on $X$ for which $\int \phi(x) \ d\mu(x) = 1$. Then there is an
  $E$-invariant $\mu$-conull \Borel set $C \subseteq X$ for which there is
  a \Borel retraction $\psi \from C \to C \intersection \support{\phi}$ such that
  $\graph{\psi} \subseteq E$ and $\cardinality{\preimage{\psi}{x}} = \phi(x)$
  for all $x \in C$.
\end{proposition}
  
\begin{propositionproof}
  As the fact that $\int \phi(x) \ d\mu(x) < \infty$ ensures that $\preimage
  {\phi}{\infty}$ is $\mu$-null, the $E$-quasi-invariance of $\mu$ implies
  that so too is its $E$-saturation. By throwing out this set, we can assume
  that $\image{\phi}{X} \subseteq \N$.
    
  By \cite[Lemma 7.3]{KechrisMiller}, there is a maximal \Borel set $\calF$ of pairwise
  disjoint finite subsets $F$ of $E$-classes such that
  $\cardinality{F \intersection \support{\phi}} = 1$ and $\cardinality{F} = 
  \phi(x)$, where $x \in F \intersection 
  \support{\phi}$. The uniformization theorem for \Borel subsets of the
  plane with countable vertical sections ensures that the set $D = \union
  [\calF]$ is \Borel, as is the function $\psi \from D \to D$ given by $\psi(x)
  = y \iff \exists F \in \calF \ (x \in F \mathand y \in F \intersection \support
  {\phi})$.
  
  It only remains to show that there is an $E$-invariant $\mu$-conull
  \Borel set $C \subseteq D$. As $\mu$ is $E$-quasi-invariant, it is
  sufficient to show that $D$ is $\mu$-conull. Suppose, towards a contradiction,
  that $X \setminus D$ is $\mu$-positive. As $\mu$ is $E$-invariant,
  the uniformization theorem for \Borel subsets of the plane with countable
  vertical sections ensures that $\int_D \phi(x) \ d\mu(x) = \mu(D)$. The
  fact that $\int \phi(x) \ d\mu(x) = 1$ therefore implies
  that $\int_{\support{\phi}\setminus D} \phi(x) \ d\mu(x) = \int_{X \setminus
  D} \phi(x) \ d\mu(x) = \mu(X \setminus D)$. In particular, it follows that
  $\support{\phi} \setminus D$ is also $\mu$-positive. As the uniformization
  theorem for \Borel subsets of the plane with countable vertical sections ensures
  that the $E$-saturation of this set is \Borel, the $E$-ergodicity of
  $\mu$ implies that this $E$-saturation is $\mu$-conull. And since the maximality
  of $\calF$ ensures that $\cardinality{((X \setminus \support{\phi}) \setminus D)
  \intersection \equivalenceclass{x}{E}} < \phi(x) - 1$ for all $x \in \support{\phi}
  \setminus D$, the $E$-invariance of $\mu$ along with the uniformization theorem
  for \Borel subsets of the plane with countable vertical sections implies that $\int_{X
  \setminus D} \phi(x) \ d\mu(x) > \mu(X \setminus D)$, the desired contradiction.
\end{propositionproof}

We next examine when binary relations can be shifted along equivalence
relations so as to stabilize the cardinalities of their sections.

For each $x \in X$, the \definedterm{\textrm{$x$}{th} vertical section} of a set
$R \subseteq X \times Y$ is given by $\verticalsection{R}{x} = \set{y \in Y}
[\pair{x}{y} \in R]$.

\begin{proposition} \label{ergodicgenerator:matching:weak}
  Suppose that $X$ and $Y$ are \Polish spaces, $E$ is a countable
  \Borel equivalence relation on $X$, $R \subseteq X \times Y$ is a \Borel set
  with countable vertical sections, and $\mu$ is an $E$-ergodic $E$-invariant
  \Borel probability measure on $X$ for which $\int \cardinality{\verticalsection
  {R}{x}} \ d\mu(x) = 1$. Then there is a \Borel function $\pi \from R \to X$ such
  that $\forall \pair{x}{y} \in R \ x \mathrel{E} \pi(x, y)$ and $\cardinality
  {\verticalsection{\image{(\pi \times \projection{Y})}{R}}{x}} = 1$ $\mu$-almost everywhere.
\end{proposition}

\begin{propositionproof}
  The uniformization theorem for \Borel subsets of the plane with countable vertical
  sections ensures that the function $\phi \from X \to \N \union \set{\infty}$ given
  by $\phi(x) = \cardinality{\verticalsection{R}{x}}$ is \Borel. Proposition \ref
  {ergodicgenerator:manytoone} therefore yields an $E$-invariant $\mu$-conull
  \Borel set $C \subseteq X$ for which there is a \Borel retraction $\psi \from C \to
  C \intersection \image{\projection{X}}{R}$, whose graph is contained in $E$, such that
  $\cardinality{\preimage{\psi}{x}} = \cardinality{\verticalsection{R}{x}}$ for all $x
  \in C$.
  
  By the uniformization theorem for \Borel subsets of the plane with countable
  vertical sections, there are \Borel functions $\rho_n \from X \to X$ and
  $\psi_n \from C \to C$ such that $\sequence{\rho_n(x)}
  [n < \cardinality{\verticalsection{R}{x}}]$ is an injective enumeration of
  $\verticalsection{R}{x}$ for all $x \in X$, and $\sequence{\psi_n(x)}[n <
  \cardinality{\preimage{\psi}{x}}]$ is an injective enumeration of
  $\preimage{\psi}{x}$ for all $x \in C$. Then the function $\pi \from R \to X$
  given by
  \begin{equation*}
    \pi(x, \rho_n(x)) =
      \begin{cases}
        \psi_n(x) & \text{if $x \in C$ and $n < \cardinality{\verticalsection{R}{x}}$,
          and} \\
        x & \text{otherwise}
      \end{cases}
  \end{equation*}
  is as desired.
\end{propositionproof}

For each $y \in Y$, the \definedterm{\textrm{$y$}{th} horizontal section} of a
set $R \subseteq X \times Y$ is given by $\horizontalsection{R}{y} = \set{x \in X}
[\pair{x}{y} \in R]$.

\begin{proposition} \label{ergodicgenerator:matching}
  Suppose that $X$ is a \Polish space, $E$ is a countable \Borel equivalence
  relation on $X$, $R \subseteq X \times X$ is a \Borel set with countable
  horizontal and vertical sections, and $\mu$ is an $E$-ergodic $E$-invariant
  \Borel probability measure on $X$ with the property that $\int \cardinality{\verticalsection
  {R}{x}} \ d\mu(x) = \int \cardinality{\horizontalsection{R}{x}} \ d\mu(x) = 1$.
  Then there are \Borel functions $\pi_X, \pi_Y \from R \to X$ such that
  $\forall \pair{x}{y} \in R \ (x \mathrel{E} \pi_X(x, y) \mathand y \mathrel{E} \pi_Y(x, y))$
  and $\cardinality{\verticalsection{\image{(\pi_X \times \pi_Y)}{R}}{x}} = \cardinality
  {\horizontalsection{\image{(\pi_X \times \pi_Y)}{R}}{x}} = 1$ $\mu$-almost everywhere.
\end{proposition}

\begin{propositionproof}
  This follows from two applications of Proposition \ref{ergodicgenerator:matching:weak}.
\end{propositionproof}

We next observe that if an invariant \Borel probability measure is ergodic
with respect to a subequivalence relation generated by a free \Borel action of a countable
non-abelian free group, then by passing to an appropriate subequivalence relation, we can
assume that it is ergodic with respect to the subequivalence relation generated by one of
the generators.

We say that an equivalence relation is \definedterm{aperiodic} if all of its classes are
infinite. If $E$ is an aperiodic countable \Borel equivalence relation and $\mu$ is an
$E$-invariant \Borel probability measure, then the uniformization theorem for \Borel
subsets of the plane with countable vertical sections ensures that $\mu(A) = 0$ for
every \Borel set $A \subseteq X$ intersecting each $E$-class in finitely-many points.

\begin{proposition} \label{ergodicgenerator:onemeasure}
  Suppose that $X$ is a \Polish space, $E$ is a countable \Borel equivalence relation on
  $X$, $\mu$ is an $E$-invariant \Borel probability measure on $X$, and there is a free
  \Borel action of a countable non-abelian free group on $X$ generating a subequivalence
  relation of $E$ with respect to which $\mu$ is ergodic. Then for every non-abelian group
  $\Gamma$ freely generated by a countable set $S$ and for every $\gamma \in S$, there
  is a free \Borel action of $\Gamma$ on $X$ generating a subequivalence relation of $E$
  such that $\mu$ is ergodic with respect to the equivalence relation generated by $\gamma$.
\end{proposition}

\begin{propositionproof}
  Note that if $\pair{\gamma_0}{\gamma_1}$ freely generates $\F{2}$, then $\sequence
  {\gamma_1^n \gamma_0 \gamma_1^{-n}}[n \in \N]$ freely generates a copy of $\F
  {\aleph_0}$ within $\F{2}$. Thus we need only construct the desired action for $\F{2}$,
  and can therefore freely discard $E$-invariant $\mu$-null \Borel sets in the course of
  the construction.

  A \definedterm{graph} on $X$ is an irreflexive symmetric subset $G$ of $X \times X$.
  The equivalence relation \definedterm{generated} by such a graph is the smallest
  equivalence relation on $X$ containing $G$, and a \definedterm{graphing} of $E$ is a
  graph generating $E$.

  The \definedterm{$\mu$-cost} of a locally countable \Borel graph $G$ on $X$ is given by
  $\cost{\mu}{G} = \frac{1}{2} \int \cardinality{\verticalsection{G}{x}} \ d\mu(x)$, and the 
  \definedterm{$\mu$-cost} of $E$ is given by
  \begin{equation*}
    \cost{\mu}{E} = \inf \set{\cost{\mu}{G}}[G \text{ is a \Borel graphing of } E].
  \end{equation*}
  
  As we can assume that $E$ is itself generated by a free \Borel action of a countable
  non-abelian free group, \Gaboriau's formula for the cost of such a relation ensures that
  $\cost{\mu}{E} \ge 2$ (see, for example, \cite[Theorem 27.6]{KechrisMiller}).
  
  Fix an aperiodic hyperfinite \Borel subequivalence relation $E_0$ of $E$ with respect to
  which $\mu$ is ergodic (see, for example, \cite[Lemma 9.3.2]{Zimmer}). Then there is a
  \Borel automorphism $T_0 \from X \to X$ generating $E_0$ (see, for example, \cite
  [Theorem 5.1]{DoughertyJacksonKechris}).

  As equality on $X$ is \Borel, we obtain a \Borel graphing $G$ of $E$ by subtracting this
  relation from $E$. A \definedterm{path} through a graph $H$ is a sequence of the form
  $\sequence{x_i}[i \le n]$, where $n \in \N$ and $\pair{x_i}{x_{i+1}} \in H$ for all $i < n$. Such
  a path is a \definedterm{cycle} if $n \ge 3$, $\sequence{x_i}[i < n]$ is injective, and $x_0 =
  x_n$. A graph is \definedterm{acyclic} if it has no cycles. The cycle-cutting lemma of
  \Kechris-\Miller and \Pichot (see, for example, \cite[Lemma 28.11]{KechrisMiller}) yields
  an acyclic \Borel subgraph $H$ of $G$ with the property that $\graph{T_0} \subseteq H$
  and $\cost{\mu}{H} \ge \cost{\mu}{E}$, thus $\cost{\mu}{H \setminus \graph{T_0^{\pm 1}}} \ge 1$.
  
  An \definedterm{oriented graph} on $X$ is an antisymmetric irreflexive subset $K$ of $X
  \times X$. The graph \definedterm{induced} by such an oriented graph is the smallest graph
  containing $K$. The fact that $X$ admits a \Borel linear order easily yields a \Borel oriented
  graph $K$ on $X$ generating $H \setminus \graph{T_0^{\pm 1}}$.
 
  The \definedterm{in-degree} of $K$ at a point $y$ is given by $\cardinality{\horizontalsection
  {K}{y}}$, whereas the \definedterm{out-degree} of $K$ at a point $x$ is given by
  $\cardinality{\verticalsection{K}{x}}$. As the $E$-invariance of $\mu$ ensures that the average
  in-degree of $K$ with respect to $\mu$ equals the average out-degree of $K$ with respect to
  $\mu$, they must be at least one. As $\mu$ is necessarily continuous, the uniformization
  theorem for \Borel subsets of the plane with countable vertical sections ensures the existence
  of a \Borel oriented subgraph $L$ of $K$ with the property that the average in-degree of $L$
  with respect to $\mu$ and the average out-degree of $L$ with respect to $\mu$ are
  exactly one.
  
  By Proposition \ref{ergodicgenerator:matching}, there are \Borel functions $\pi_X, \pi_Y
  \from L \to X$ such that $\forall \pair{x}{y} \in L \ (x \mathrel{E} \pi_X(x, y) \mathand
  y \mathrel{E} \pi_Y(x, y))$ and the in-degree and out-degree of $\image{(\pi_X
  \times \pi_Y)}{L}$ at $\mu$-almost every point is one. As $\mu$ is $E$-quasi-invariant,
  the uniformization theorem for \Borel subsets of the plane with countable vertical
  sections ensures that, after throwing out an $E$-invariant $\mu$-null \Borel set,
  we can assume that the in-degree and out-degree of $\image{(\pi_X \times \pi_Y)}
  {L}$ at every point is one. Let $T_1$ be the unique \Borel automorphism whose
  graph is $\image{(\pi_X \times \pi_Y)}{L}$. As $H$ is acyclic, the requirement that
  $\forall \pair{x}{y} \in L \ (x \mathrel{E} \pi_X(x, y) \mathand y \mathrel{E} \pi_Y(x, y))$
  ensures that the action of $\F{2}$ generated by $T_0$ and $T_1$ is free.
\end{propositionproof}

\begin{remark} \label{ergodicgenerator:onemeasure:remark}
  After throwing out an $E$-invariant $\mu$-null \Borel set, the conclusion of Proposition
  \ref{ergodicgenerator:onemeasure} can be established from the weaker hypothesis that
  there is a \Borel subequivalence relation $F$ of $E$, for which $\cost{\mu}{F} > 1$, with 
  respect to which $\mu$ is ergodic. To see this, note that we can assume that $\cost{\mu}
  {E} > 1$. As $\mu$ is necessarily continuous, there is a \Borel set $B \subseteq X$ for
  which there is a natural number $k \ge 1 / (\cost{\mu}{E} - 1)$ such that $\mu(B) = 1 / k$,
  thus $\mu(B) \le \cost{\mu}{E} - 1$. \Gaboriau's formula for the cost of the restriction of a
  countable \Borel equivalence relation (see, for example, \cite[Theorem 21.1]
  {KechrisMiller}) then ensures that $\cost{\restriction{\mu}{B}}{\restriction{E}{B}} =
  (\cost{\mu}{E} - 1) + \mu(B) \ge 2\mu(B)$. The proof of Proposition \ref
  {ergodicgenerator:onemeasure} therefore yields the desired action of $\F{2}$, albeit on
  $B$, not on $X$. To rectify this problem, note that by the proof of \cite[Lemma 7.10]
  {KechrisMiller}, after throwing out an $E$-invariant $\mu$-null \Borel set, we can
  assume that there is a partition of $X$ into \Borel sets $B_1, \ldots, B_k$, as well as
  \Borel isomorphisms $\pi_i \from B \to B_i$ whose graphs are contained in $E$, for
  all $1 \le i \le k$. Let $T_0'$ denote the \Borel automorphism of $X$ which agrees
  with $\pi_{i+1} \composition \inverse{\pi_i}$ on $B_i$, for all $1 \le i < k$, and which
  agrees with $\pi_1 \composition T_0 \composition \inverse{\pi_k}$ on $B_k$. Let
  $T_1'$ denote the \Borel automorphism of $X$ which agrees with $\pi_i \composition
  T_1^i \composition T_0 \composition T_1^{-i} \composition \inverse{\pi_i}$ on
  $B_i$, for all $1 \le i \le k$. Then $T_0'$ and $T_1'$ yield the desired action of $\F{2}$.
\end{remark}

We equip the space of \Borel probability measures $\mu$ on a \Polish space $X$
with the (standard) \Borel structure generated by the functions $\evaluation_B(\mu) = \mu(B)$,
where $B$ varies over all \Borel subsets of $X$.

A \definedterm{uniform ergodic decomposition} of a countable \Borel equivalence relation
$E$ on a \Polish space $X$ is a sequence $\sequence{\mu_x}[x \in X]$ of
$E$-ergodic $E$-invariant \Borel probability measures on $X$ such that (1) $\mu_x = \mu_y$
whenever $x \mathrel{E} y$, (2) $\mu(\set{x \in X}[\mu = \mu_x]) = 1$ for every
$E$-ergodic $E$-invariant \Borel probability measure $\mu$ on $X$, and (3) $\mu = \int \mu_x
\ d\mu(x)$ for every $E$-invariant \Borel probability measure $\mu$ on $X$. The
\Farrell-\Varadarajan uniform ergodic decomposition theorem (see, for example, \cite[Theorem
3.3]{KechrisMiller}) ensures the existence of \Borel such decompositions. We next establish
an ergodicity-free analog of Proposition \ref{ergodicgenerator:onemeasure}, by uniformly
pasting together actions obtained from the latter along such a decomposition.

\begin{proposition} \label{ergodicgenerator:manymeasures}
  Suppose that $X$ is a \Polish space, $E$ is a countable \Borel equivalence relation
  on $X$, $\mu$ is an $E$-invariant \Borel probability measure on $X$, $\sequence
  {\mu_x}[x \in X]$ is a \Borel uniform ergodic decomposition of $E$, and for
  $\mu$-almost all $x \in X$ there is a free \Borel action of a countable non-abelian free
  group generating a subequivalence relation of $E$ with respect to which $\mu_x$ is
  ergodic. Then for every non-abelian group $\Gamma$ freely generated by a countable
  set $S$ and for every $\gamma \in S$, there is a free \Borel action of $\Gamma$ on $X$
  generating a subequivalence relation of $E$ such that for $\mu$-almost all $x \in X$,
  the measure $\mu_x$ is ergodic with respect to the equivalence relation generated by
  $\gamma$.
\end{proposition}

\begin{propositionproof}
  As before, it is sufficient to construct the desired action for $\F{2}$, freely discarding
  $E$-invariant $\mu$-null \Borel sets as we proceed. Fix $\pair{\gamma_0}{\gamma_1}$
  freely generating $\F{2}$.
  
  Fix a countable basis $\set{U_n}[n \in \N]$ for $X$ which is closed under finite
  unions, and define $G \subseteq \functions{\N \times \N}{\N} \times X$ by
  $\verticalsection{G}{\phi} = \intersection[m \in \N][{\union[n \in \N][U_{\phi(m, n)}]}]$.
  By the uniformization theorem for \Borel subsets of the plane with countable vertical
  sections, there are \Borel functions $f_n \from X \to X$ such that $E = \union[n \in \N]
  [\graph{f_n}]$.
  
  Define $R \subseteq \functions{\N}{(\functions{\N \times \N}{\N})} \times (X \times X)$
  by $\verticalsection{R}{\phi} = \union[n \in \N][\graph{\restriction{f_n}{\verticalsection
  {G}{\phi(n)}}}]$ for all $\phi \in \functions{\N}{(\functions{\N \times \N}{\N})}$, and let
  $B$ denote the \Borel set consisting of all $\pair{\pair{\phi_0}{\phi_1}}{x} \in
  (\functions{\N}{(\functions{\N \times \N}{\N})} \times \functions{\N}{(\functions{\N \times
  \N}{\N})}) \times X$ with the property that the sets $\verticalsection{R}{\phi_0} \intersection
  (\equivalenceclass{x}{E} \times \equivalenceclass{x}{E})$ and $\verticalsection{R}
  {\phi_1} \intersection (\equivalenceclass{x}{E} \times \equivalenceclass{x}{E})$ are
  graphs of functions inducing a free action of $\F{2}$ on $\equivalenceclass{x}{E}$. For
  each pair $\phi = \pair{\phi_0}{\phi_1}$ in $\functions{\N}{(\functions{\N \times \N}{\N})}
  \times \functions{\N}{(\functions{\N \times \N}{\N})}$, let $\cdot_\phi$ denote the action
  of $\F{2}$ on $\verticalsection{B}{\phi}$ for which the graph of the function $x \mapsto
  \gamma_i \cdot_\phi x$ is $\verticalsection{R}{\phi_i} \intersection (\verticalsection{B}
  {\phi} \times \verticalsection{B}{\phi})$, for all $i < 2$. Let $F$ denote the equivalence
  relation on $B$ generated by the function $\pair{\phi}{x} \mapsto \pair{\phi}{\gamma_0
  \cdot_\phi x}$, and fix a \Borel uniform ergodic decomposition $\sequence{\nu_{\phi,x}}
  [\pair{\phi}{x} \in B]$ of $F$.
  
  A subset of a \Polish space is \definedterm{analytic} if it is the image of a
  \Borel subset of a \Polish space under a \Borel function. We use $\sigmaclass
  {\Sigmaclass[1][1]}$ to denote the smallest $\sigma$-algebra containing all such sets,
  and we say that a function is \definedterm{$\sigmaclass{\Sigmaclass[1][1]}$-measurable}
  if pre-images of open sets are $\Sigmaclass[1][1]$. Define $\projection{X} \from X \times Y
  \to X$ by $\projection{X}(x, y) = x$. A \definedterm{uniformization} of a
  set $R \subseteq X \times Y$ is a function $\phi \from \image{\projection{X}}{R} \to Y$
  whose graph is contained in $R$.
  
  Let $\Lambda$ denote the push-forward of $\mu$ through the function $x \mapsto
  \mu_x$. The regularity of \Borel probability measures on \Polish spaces ensures that if
  $\nu$ is an $E$-invariant \Borel probability measure on $X$ concentrating on
  $\verticalsection{B}{\phi}$, then every free \Borel action of $\F{2}$ on an $E$-invariant
  $\nu$-conull \Borel subset of $X$ agrees with some $\cdot_\phi$ on an $E$-invariant
  $\nu$-conull \Borel subset of $\verticalsection{B}{\phi}$. And $\restriction{\nu}
  {\verticalsection{B}{\phi}}$ is ergodic with respect to the equivalence relation generated
  by the action of $\gamma_0$ on $\verticalsection{B}{\phi}$ if and only if the measures
  induced by $\nu$ and $\restriction{\nu_{\phi,x}}{(\singleton{\phi}\times \verticalsection
  {B}{\phi}})$ on $\verticalsection{B}{\phi}$ agree, for $\nu$-almost all $x \in X$. As the
  set of pairs $\pair{\nu}{\phi}$ satisfying this latter property is \Borel (see, for example,
  \cite[Proposition 12.4 and Theorem 17.25]{Kechris}), the \Jankov-\vonNeumann 
  uniformization theorem for analytic subsets of the plane (see, for example, \cite[Theorem
  18.1]{Kechris}) yields a $\sigmaclass{\Sigmaclass[1][1]}$-measurable uniformization $\phi$.
  As $\domain{\phi}$ is analytic, and a result of \Lusin's ensures that every analytic set is
  $\Lambda$-measurable (see, for example, \cite[Theorem 21.10]{Kechris}), Proposition 
  \ref{ergodicgenerator:onemeasure} implies that $\domain{\phi}$ is $\Lambda$-conull. As
  $\phi$ is $\Lambda$-measurable, there is a $\Lambda$-conull \Borel set $M \subseteq
  \domain{\phi}$ on which it is \Borel.
  
  Define $C = \union[\nu \in M][{\set{x \in B_{\phi(\nu)}}[\mu_x = \nu]}]$, and 
  observe that the action of $\F{2}$ on $C$ given by $\gamma_i \cdot x = \gamma_i
  \cdot_{\phi(\mu_x)} x$, for $i < 2$, is as desired.
\end{propositionproof}

\begin{remark} \label{ergodicgenerator:manymeasures:remark}
  By Remark \ref{ergodicgenerator:onemeasure:remark}, after throwing out an
  $E$-invariant $\mu$-null \Borel set, the conclusion of Proposition \ref
  {ergodicgenerator:manymeasures} follows from the weaker assumption
  that for $\mu$-almost all $x \in X$ there is a \Borel subequivalence relation $F$ of $E$,
  for which $\cost{\mu_x}{F} > 1$, with respect to which $\mu_x$ is ergodic.
\end{remark}

\section{Stratification} \label{stratification}

Here we consider various sorts of stratifications of countable \Borel equivalence relations.
We begin by noting that strong proper inclusion of equivalence relations often gives rise
to a measure-theoretic analog.

\begin{proposition} \label{stratification:properinclusion}
  Suppose that $X$ is a \Polish space, $E$ and $F$ are countable \Borel equivalence
  relations on $X$ such that every $E$-class is properly contained in an $F$-class, and
  $\mu$ is an $E$-ergodic $F$-quasi-invariant \Borel probability measure on $X$. Then
  for every \Borel set $B \subseteq X$, the $(\restriction{E}{B})$-class of $(\restriction{\mu}
  {B})$-almost every point of $B$ is properly contained in an $(\restriction{F}{B})$-class.
\end{proposition}

\begin{propositionproof}
  The uniformization theorem for \Borel subsets of the plane with countable vertical
  sections ensures that $\saturation{B}{F}$ is \Borel, thus so too is the set $A = \set{x \in
  \saturation{B}{F}}[B \intersection \equivalenceclass{x}{E} = B \intersection
  \equivalenceclass{x}{F}]$. As the fact that every $E$-class is properly contained in an
  $F$-class ensures that $A$ is contained in the $F$-saturation of its complement, the
  $F$-quasi-invariance of $\mu$ implies that $A$ is not $\mu$-conull. As $A$ is
  $E$-invariant, the $E$-ergodicity of $\mu$ therefore ensures that $A$ is $\mu$-null. It
  only remains to note that if $x \in B \setminus A$, then $x$ is necessarily $(F \setminus
  E)$-related to some other point in $B$, thus the $(\restriction{E}{B})$-class of $x$ is
  properly contained in an $(\restriction{F}{B})$-class.
\end{propositionproof}

We say that a sequence $\sequence{E_r}[r \in \R]$ of subequivalence relations of an
equivalence relation $E$ on $X$ is a \definedterm{stratification} of $E$ if for all real
numbers $r < s$, every $E_r$-class is properly contained in an $E_s$-class. We say that
a sequence $\sequence{\cdot_r}[r \in \R]$ of actions of a group $\Gamma$ on $X$ is a 
\definedterm{stratification} of $E$ if the corresponding sequence $\sequence{E_r}[r \in
\R]$ of equivalence relations generated by the actions is a stratification.

We say that a stratification $\sequence{E_r}[r \in \R]$ of $E$ is a \definedterm
{$\mu$-stratification} if for all $\mu$-positive \Borel sets $B \subseteq X$ and all real numbers
$r < s$, the $(\restriction{E_r}{B})$-class of $(\restriction{\mu}{B})$-almost every point of $B$ is 
properly contained in an $(\restriction{E_s}{B})$-class. We say that a stratification 
$\sequence{\cdot_r}[r \in \R]$ of $E$ is a \definedterm{$\mu$-stratification} if
the corresponding sequence $\sequence{E_r}[r \in \R]$ of equivalence relations generated
by the actions is a $\mu$-stratification.

We say that a ($\mu$-)stratification $\sequence{E_r}[r \in \R]$ by equivalence relations is 
\definedterm{\Borel} if it is \Borel as a sequence of subsets of $X \times X$. Analogously, a
($\mu$-)stratification $\sequence{\cdot_r}[r \in \R]$ by actions of $\Gamma$ is \definedterm
{\Borel} if the corresponding sequence $\sequence{\graph{\cdot_r}}[r \in \R]$ is \Borel as a
sequence of subsets of $(\Gamma \times X) \times X$.

\begin{proposition} \label{stratification:ergodicinvariantactions}
  Suppose that $X$ is a \Polish space, $E$ is a countable \Borel equivalence relation on
  $X$, $\mu$ is an $E$-invariant \Borel probability measure on $X$, and
  there is a free \Borel action of a countable non-abelian free group on $X$ generating a
  subequivalence relation of $E$ with respect to which $\mu$ is ergodic. Then for every
  countable non-abelian free group $\Gamma$, there is a \Borel $\mu$-stratification of 
  $E$ by free actions of $\Gamma$ on $X$, generating equivalence relations with respect
  to which $\mu$ is ergodic.
\end{proposition}

\begin{propositionproof}
  Fix a countable set $S$ freely generating $\Gamma$, as well as
  some $\gamma \in S$. By Proposition \ref{ergodicgenerator:onemeasure}, there is a free
  \Borel action $\cdot$ of $\Gamma$ on $X$, generating a subequivalence relation of $E$, 
  such that $\mu$ is ergodic with respect to the equivalence relation generated by $\gamma$.
  Fix $\delta \in S \setminus \singleton{\gamma}$, appeal to \cite[Proposition 5.2]
  {Miller:Incomparable} to obtain a \Borel stratification $\sequence{*_r}[r \in \R]$ by actions of
  $\Z$ of the equivalence relation generated by $\delta$, and for each $r \in \R$, let $\cdot_r$
  denote the action of $\Gamma$ on $X$ given by $\delta \cdot_r x = \delta *_r x$ and
  $\lambda \cdot_r x = \lambda \cdot x$, for $\lambda \in S \setminus \singleton{\delta}$.
  Proposition \ref{stratification:properinclusion} then ensures that $\sequence{\cdot_r}[r \in \R]$
  is the desired $\mu$-stratification of $E$.
\end{propositionproof}

\begin{remark} \label{stratification:ergodicinvariantactions:remark}
  By Remark \ref{ergodicgenerator:onemeasure:remark}, after throwing out an
  $E$-invariant $\mu$-null \Borel set, the conclusion of Proposition \ref
  {stratification:ergodicinvariantactions} follows from the weaker assumption
  that there is a \Borel subequivalence relation $F$ of $E$, for which $\cost{\mu}{F} > 1$,
  with respect to which $\mu$ is ergodic.
\end{remark}

We next establish an analogous result in the absence of ergodicity.

\begin{proposition} \label{stratification:invariantactions}
  Suppose that $X$ is a \Polish space, $E$ is a countable \Borel equivalence relation on
  $X$, $\mu$ is an $E$-invariant \Borel probability measure on $X$, and there is a free
  \Borel action of a countable non-abelian free group on $X$ generating a subequivalence
  relation of $E$. Then for every countable non-abelian free group $\Gamma$, there is a
  \Borel $\mu$-stratification of $E$ by free actions of $\Gamma$ on $X$.
\end{proposition}

\begin{propositionproof}
  Fix a countable set $S$ freely generating a non-abelian group $\Gamma$, as well as
  some $\gamma \in S$. Clearly we can assume that $E$ is itself generated by a free
  \Borel action of a countable non-abelian free group on $X$. By the uniform ergodic
  decomposition theorem, there is a \Borel uniform ergodic decomposition $\sequence
  {\mu_x}[x \in X]$ of $E$. By Proposition \ref{ergodicgenerator:manymeasures}, there
  is a free \Borel action $\cdot$ of $\Gamma$ on $X$, generating a subequivalence relation of
  $E$, such that for $\mu$-almost all $x \in X$, the measure $\mu_x$ is ergodic with
  respect to the equivalence relation generated by $\gamma$. Fix $\delta \in S \setminus 
  \singleton{\gamma}$, appeal to \cite[Proposition 5.2]{Miller:Incomparable} to obtain a
  \Borel stratification $\sequence{*_r}[r \in \R]$ by actions of $\Z$ of the equivalence relation
  generated by $\delta$, and for each $r \in \R$, let $\cdot_r$ denote the action of $\Gamma$ on
  $X$ given by $\delta \cdot_r x = \delta *_r x$ and $\lambda \cdot_r x = \lambda \cdot x$, for
  $\lambda \in S \setminus \singleton{\delta}$. Proposition \ref{stratification:properinclusion}
  ensures that $\sequence{\cdot_r}[r \in \R]$ is a $\mu_x$-stratification of $E$ for $\mu$-almost
  all $x \in X$, and is therefore a $\mu$-stratification of $E$.
\end{propositionproof}

\begin{remark} \label{stratification:invariantactions:remark}
  After throwing out an $E$-invariant $\mu$-null \Borel set, the conclusion of Proposition 
  \ref{stratification:invariantactions} follows from the weaker assumption that for all $\mu$-positive
  \Borel sets $B \subseteq X$, there is a \Borel subequivalence relation $F$ of $\restriction{E}{B}$
  for which $\cost{\restriction{\mu}{B}}{F} > \mu(B)$. In light of Remark \ref
  {ergodicgenerator:manymeasures} and the proof of Proposition \ref{stratification:invariantactions},
  to see this, it is sufficient to show that there is a \Borel subequivalence relation $F$ of $E$ with
  the property that if $\sequence{\mu_x}[x \in X]$ is an ergodic decomposition of $F$, then
  $\cost{\mu_x}{F} > 1$ for $\mu$-almost all $x \in X$. In fact, it is enough to produce such an $F$ on 
  an $E$-invariant $\mu$-positive \Borel subset of $X$. Towards this end, fix a \Borel subequivalence 
  relation $F'$ of $E$ such that $\cost{\mu}{F'} > 1$, and appeal to the uniform ergodic decomposition 
  theorem to produce a \Borel uniform ergodic decomposition $\sequence{\mu_x}[x \in X]$ of $F'$. As
  \cite[Proposition 18.1]{KechrisMiller} ensures that the set $C = \set{x \in X}[\cost{\mu_x}{F'}
  > 1]$ is co-analytic, and every co-analytic subset of $X$ is $\mu$-measurable, the cost integration 
  formula (see \cite[Corollary 18.6]{KechrisMiller}) ensures that it is $\mu$-positive. Fix a
  $\mu$-positive \Borel set $B \subseteq X$ contained in $C$. As $\mu$ is $F'$-invariant, we can 
  assume that $B$ is $F'$-invariant. By the uniformization theorem for \Borel subsets of the plane
  with countable vertical sections, the set $\saturation{B}{E}$ is \Borel and there is a \Borel
  function $\phi \from \saturation{B}{E} \to B$ whose graph is contained in $E$. Then the formula
  for the cost of a restriction of a countable \Borel equivalence relation ensures that the
  equivalence relation $F$ on $\saturation{B}{E}$ given by $x \mathrel{F} y \iff \phi(x) \mathrel{F'} 
  \phi(y)$ is as desired. 
\end{remark}

A \definedterm{treeing} of an equivalence relation is an acyclic graphing. We
say that a \Borel equivalence relation $E$ on $X$ is \definedterm{treeable} if it
has a \Borel treeing. We say that a countable \Borel equivalence relation $E$
on $X$ is \definedterm{compressible} if there is a \Borel injection $\phi \from X
\to X$, whose graph is contained in $E$, such that $X \setminus \image{\phi}{X}$
intersects every $E$-class.

\begin{proposition} \label{stratification:compressiblefreeaction}
  Suppose that $X$ is a \Polish space, $E$ is a compressible treeable countable \Borel 
  equivalence relation on $X$, and $\Gamma$ is a countable non-abelian free group. Then
  there is a free \Borel action of $\Gamma$ on $X$ generating $E$.
\end{proposition}

\begin{propositionproof}
  A straightforward modification of the proof of \cite[Corollary 3.11]{JacksonKechrisLouveau}
  reveals that there is a \Borel treeing $G$ of $E$ which is generated by a \Borel automorphism
  $T_0 \from X \to X$ and a \Borel isomorphism $T_1 \from A \to B$, where $A, B \subseteq X$ are
  disjoint \Borel sets. Set $C = X \setminus (A \union B)$, and for each $x \in X$, let $D(x)$ denote
  the unique set in $\set{A, B, C}$ containing $x$.
  
  Suppose now that $k \in \set{2, 3, \ldots, \aleph_0}$ and $\Gamma$ is freely generated by
  $\sequence{\gamma_i}[i < k]$. The idea is to again use compressibility to replace $E$ with
  $E \times \square{\N}$, and to use the right-hand coordinate to accomodate the generators
  other than the first two, as well as the points at which $T_1$ and $\inverse{T_1}$ are not
  defined.
  
   Towards this end, fix bijections $\gamma_0^D \from \N \setminus \singleton{0} \to 
  \N \setminus \singleton{0}$ for all $D \in \set
  {A, B, C}$, $\gamma_1^A \from \N \setminus \singleton{0} \to \N$, $\gamma_1^B \from \N \to
  \N \setminus \singleton{0}$, $\gamma_1^C \from \N \to \N$, as well as $\gamma_i^D \from \N \to \N$
  for all $D \in \set{A, B, C}$ and $1 < i < k$, with the property that for all $D \in \set{A, B, C}$, the corresponding 
  approximation to an action of $\F{k}$ on $\N$, given by
  \begin{equation*}
    (\gamma_{s(0)}^{t(0)} \cdots \gamma_{s(n)}^{t(n)})^D \cdot x = ((\gamma_{s(0)}^D)^{t(0)} 
      \composition \cdots \composition (\gamma_{s(n)}^D)^{t(n)})(x),
  \end{equation*}
  is both \definedterm{free} and \definedterm{transitive}, in the sense that:
  \begin{enumerate}
    \item $\forall n \in \N \forall \gamma \in \F{k} \ (\gamma^D \cdot n = n \iff \gamma = \id)$.
    \item $\forall m, n \in \N \exists \gamma \in \F{k} \ (\gamma^D \cdot m = n)$.
  \end{enumerate}
  Let $\gamma_0$ act on $X \times \N$ via
  \begin{equation*}
    \gamma_0 \cdot \pair{x}{n} =
      \begin{cases}
        \pair{T_0(x)}{n} & \text{if $n = 0$, and} \\
        \pair{x}{\gamma_0^{D(x)} \cdot n} & \text{otherwise.}
      \end{cases}
  \end{equation*}
  Similarly, let $\gamma_1$ act on $X \times \N$ via
  \begin{equation*}
    \gamma_1 \cdot \pair{x}{n} =
      \begin{cases}
        \pair{T_1(x)}{n} & \text{if $n = 0$ and $x \in A$, and} \\
        \pair{x}{\gamma_1^{D(x)} \cdot n} & \text{otherwise.}
      \end{cases}
  \end{equation*}
  And finally, let $\gamma_i$ act on $X \times \N$ via $\gamma_i \cdot \pair{x}{n} =
  \pair{x}{\gamma_i^{D(x)} \cdot n}$, for all $1 < i < k$. This defines a free \Borel action of
  $\Gamma$ generating $E \times \square{\N}$, so the proposition follows from the
  fact that the compressibility of $E$ is equivalent to the existence of a \Borel isomorphism
  between $E$ and $E \times \square{\N}$ (see, for example, \cite[Proposition 2.5]
  {DoughertyJacksonKechris}).
\end{propositionproof}

In particular, this yields the following stratification result.

\begin{proposition} \label{stratification:compressibleequivalencerelations}
  Suppose that $X$ is a \Polish space, $E$ is a compressible treeable countable \Borel
  equivalence relation on $X$, and $\mu$ is a \Borel probability measure on $X$ such that 
  $E$ is $\mu$-nowhere hyperfinite. Then there is a \Borel $\mu$-stratification of $E$ by 
  $\mu$-nowhere hyperfinite equivalence relations, each of which is generated by free 
  \Borel actions of every countable non-abelian free group on $X$.
\end{proposition}

\begin{propositionproof}
  Fix a \Borel isomorphism $\pi \from X \to X \times \N$ of $E$ with $E \times \square{\N}$,
  and let $\nu$ denote the push-forward of $\mu$ through $\pi$. Fix an $(E \times \square
  {\N})$-quasi-invariant \Borel probability measure $\nu' \gg \nu$ on $X \times \N$ (see, for
  example, the proof of \cite[Corollary 10.2]{KechrisMiller}). Define a \Borel measure $\nu_0'$
  on $X$ by $\nu_0'(B) = \nu'(B \times \singleton{0})$, for all \Borel sets $B \subseteq X$. By 
  \cite[Theorem 5.7]{ConleyMiller}, there is a \Borel $\mu$-stratification $\sequence{E_r}[r \in
  \R]$ of $E$ by $\nu_0'$-nowhere hyperfinite countable \Borel equivalence relations. Then the 
  pullback of $\sequence{E_r \times \square{\N}}[r \in \R]$ through $\pi$ yields a \Borel
  stratification of $E$ by $\mu$-nowhere hyperfinite countable \Borel equivalence relations. As
  these relations are necessarily compressible, Proposition \ref
  {stratification:compressiblefreeaction} implies that they are induced by free \Borel actions of
  every countable non-abelian free group on $X$.
\end{propositionproof}

The following ensures that \Borel $\mu$-stratifications by equivalence relations generated by
free \Borel actions of countable groups give rise to \Borel $\mu$-stratifications by actions.

\begin{proposition} \label{stratification:equivalencerelationstoactions}
  Suppose that $X$ is a \Polish space, $E$ is a countable \Borel equivalence relation on
  $X$, $\Gamma$ is a countable group, $\mu$ is a \Borel probability measure on $X$, and
  $\sequence{E_r}[r \in \R]$ is a \Borel $\mu$-stratification of $E$ by equivalence relations
  generated by free \Borel actions of $\Gamma$ on $X$. Then there is a \Borel $\mu$-stratification
  $\sequence{\cdot_r}[r \in \R]$ of $E$ by free \Borel actions of $\Gamma$ on $X$ for which
  there is a \Borel function $\pi \from \R \to \R$ such that $E_{\pi(r)}$ is the equivalence relation
  generated by $\cdot_r$ on an $E$-invariant $\mu$-conull \Borel set, for all $r \in \R$.
\end{proposition}

\begin{propositionproof}
  Fix a countable basis $\set{U_n}[n \in \N]$ for $X$ which is closed under finite unions, and
  define $G \subseteq \functions{\N \times \N}{\N} \times X$ by $\verticalsection{G}{\phi} = 
  \intersection[m \in \N][{\union[n \in \N][U_{\phi(m, n)}]}]$. By the uniformization theorem for
  \Borel subsets of the plane with countable vertical sections, there are \Borel functions $f_n
  \from X \to X$ such that $E = \union[n \in \N][\graph{f_n}]$.
  
  Define $R \subseteq \functions{\N}{(\functions{\N \times \N}{\N})} \times (X \times X)$ by
  $\verticalsection{R}{\phi} = \union[n \in \N][\graph{\restriction{f_n}{\verticalsection{G}
  {\phi(n)}}}]$, and let $B$ denote the set of all pairs $\pair{\phi}{x} \in \functions{\Gamma}
  {(\functions{\N}{(\functions{\N \times \N}{\N})})} \times X$ for which the sets of the form
  $\verticalsection{R}{\phi(\gamma)} \intersection (\equivalenceclass{x}{E} \times
  \equivalenceclass{x}{E})$ are graphs of functions inducing a free action of $\Gamma$ on
  $\equivalenceclass{x}{E}$. For each function $\phi \from \Gamma \to \functions{\N}
  {(\functions{\N \times \N}{\N})}$, let $\cdot_\phi$ denote the action of $\Gamma$ on
  $\verticalsection{B}{\phi}$ with the property that the graph of the function $x \mapsto
  \gamma \cdot_\phi x$ is $\verticalsection{R}{\phi(\gamma)} \intersection 
  (\verticalsection{B}{\phi} \times \verticalsection{B}{\phi})$, for all $\gamma \in \Gamma$.
  
  The regularity of \Borel probability measures on \Polish spaces ensures that for all $r \in \R$,
  the equivalence relation generated by an action of the form $\cdot_\phi$ is $E_r$ on an
  $E$-invariant $\mu$-conull \Borel subset of $\verticalsection{B}{\phi}$. As the set of pairs 
  $\pair{r}{\phi}$ satisfying this latter property is \Borel, the uniformization theorem for analytic
  subsets of the plane yields a $\sigmaclass{\Sigmaclass[1][1]}$-measurable uniformization
  $\phi$. Fix a continuous \Borel probability measure $m$ on $\R$. As $\phi$ is
  necessarily $m$-measurable, there is an $m$-conull \Borel set $R \subseteq \R$ on which
  $\phi$ is \Borel. As $R$ is necessarily uncountable, the proof of the perfect set theorem (see,
  for example, \cite[Theorem 13.6]{Kechris}) yields an order-preserving continuous embedding
  of $\Cantorspace$ (equipped with the lexicographic order) into $R$. And by composing such
  an embedding with any order-preserving \Borel embedding of $\R$ into $\Cantorspace$, we
  obtain an order-preserving \Borel embedding $\pi$ of $\R$ into $R$, in which case the
  sequence $\sequence{\cdot_{\pi(r)}}[r \in \R]$ is as desired.
\end{propositionproof}

Finally, we establish an analog of Proposition \ref{stratification:invariantactions} without invariance.

\begin{proposition} \label{stratification:actions}
  Suppose that $X$ is a \Polish space, $E$ is a countable \Borel equivalence relation on
  $X$, $\mu$ is a \Borel probability measure on $X$, and there is a free \Borel action of a
  countable non-abelian free group on $X$ generating a $\mu$-nowhere hyperfinite
  subequivalence relation of $E$. Then for every countable non-abelian free group $\Gamma$,
  there is a \Borel $\mu$-stratification of $E$ by free actions of $\Gamma$ on $X$ generating
  $\mu$-nowhere hyperfinite equivalence relations.
\end{proposition}

\begin{propositionproof}
  As a result of \Hopf's ensures that $X$ is compressible off of an $E$-invariant \Borel set
  on which $\mu$ is equivalent to an $E$-invariant \Borel probability measure (see, for
  example, \cite[\S10]{Nadkarni}), the desired result is a consequence of Propositions \ref
  {stratification:invariantactions}, \ref{stratification:compressibleequivalencerelations}, and
  \ref{stratification:equivalencerelationstoactions}.
\end{propositionproof}

\begin{remark} \label{stratification:actions:remark}
  By Remark \ref{stratification:invariantactions:remark}, after throwing out an $E$-invariant
  $\mu$-null \Borel set, the conclusion of Proposition \ref{stratification:actions} follows from
  the weaker assumption that for all $\mu$-positive \Borel subsets $B \subseteq X$ such
  that $\restriction{\mu}{B}$ is equivalent to an $E$-invariant \Borel probability measure $\nu$,
  there is a \Borel subequivalence relation $F$ of $\restriction{E}{B}$ for which
  $\cost{\restriction{\nu}{B}}{F} > \nu(B)$.
\end{remark}

\section{Antichains} \label{incomparable}

We are now prepared to establish our primary results. Although these can be obtained from
the arguments of \cite{ConleyMiller} (which are themselves slight variants of arguments already
appearing in \cite{Hjorth}) by substituting the stratification results of the
previous section for those of \cite{ConleyMiller}, we will nevertheless provide the full proofs,
both for the convenience of the reader and because somewhat simpler versions are
sufficient to obtain the results we consider here.

\begin{theorem} \label{antichains:ergodicinvariantactions}
  Suppose that $X$ is a \Polish space, $E$ is a projectively separable countable
  \Borel equivalence relation on $X$, $\mu$ is an $E$-invariant \Borel
  probability measure on $X$, and there is a free \Borel action of a countable non-abelian free
  group on $X$ generating a subequivalence relation of $E$ with respect to which $\mu$ is
  ergodic. Then for every countable non-abelian free group $\Gamma$, there is a \Borel
  sequence $\sequence{\cdot_r}[r \in \R]$ of free actions of $\Gamma$ on $X$, generating 
  subequivalence relations $E_r$ of $E$ with respect to which $\mu$ is ergodic, with the
  further property that $\sequence{E_r}[r \in \R]$ is an increasing sequence of relations which
  are pairwise incomparable under $\mu$-reducibility.
\end{theorem}

\begin{theoremproof}
  By Proposition \ref{stratification:ergodicinvariantactions}, there is a \Borel $\mu$-stratification
  $\sequence{\cdot_r}[r \in \R]$ of $E$ by free actions of $\Gamma$ on $X$ generating
  equivalence relations $E_r$ with respect to which $\mu$ is ergodic. By replacing
  $\sequence{\cdot_r}[r \in \R]$ with its pushforward through an order-preserving \Borel
  isomorphism of the set of all real numbers with the set of positive real numbers,
  we can assume that $\intersection[r \in \R][E_r]$ is $\mu$-nowhere hyperfinite. Let $R$ denote 
  the relation on $\R$ in which two real numbers $r$ and $s$ are related if $E_r$ is $\mu$-reducible
  to $E_s$.

  \begin{lemma} \label{antichains:ergodicinvariantactions:sections}
    Every horizontal section of $R$ is countable.
  \end{lemma}

  \begin{lemmaproof}
    Suppose, towards a contradiction, that there exists $t \in \R$ for which $\horizontalsection
    {R}{t}$ is uncountable. For each $s \in \horizontalsection{R}{t}$, fix a $\mu$-conull \Borel
    set $B_s \subseteq X$ on which there is a \Borel reduction $\phi_s \from B_s \to X$ of
    $E_s$ to $E_t$. As each $\phi_s$ is a homomorphism from $\restriction{(\intersection
    [r \in \R][E_r])}{B_s}$ to $E$, the $\mu$-nowhere hyperfiniteness of $\intersection[r \in \R]
    [E_r]$ and the projective separability of $E$ ensure the existence of distinct $r, s \in 
    \horizontalsection{R}{t}$ for which $\uniformmetric{\mu}(\phi_r, \phi_s) < 1$. Then
    $\set{x \in B_r \intersection B_s}[\phi_r(x) = \phi_s(x)]$ is a $\mu$-positive \Borel set on
    which $E_r$ and $E_s$ coincide, a contradiction.
  \end{lemmaproof}

  As $R$ is analytic (see, for example, \cite[Proposition I.15]{ConleyMiller}), a result of
  \Lusin-\Sierpinski ensures that it has the \Baire property (see, for example, \cite[Theorem
  21.6]{Kechris}). As the horizontal sections of $R$ are countable and therefore meager,
  a result of \Kuratowski-\Ulam implies that $R$ is meager (see, for example, \cite[Theorem
  8.41]{Kechris}), in which case a result of \Mycielski's yields a continuous order-preserving
  embedding $\phi \from \Cantorspace \to \R$ such that pairs of distinct sequences in
  $\Cantorspace$ are mapped to $R$-unrelated pairs of real numbers (see, for example,
  \cite[Theorem B.5]{ConleyMiller}). Fix a \Borel embedding $\psi \from \R \to \Cantorspace$
  of the usual ordering of $\R$ into the lexicographical ordering of $\Cantorspace$, and
  observe that the sequence $\sequence{\cdot_{\pi(r)}}[r \in \R]$ is as desired, where $\pi =
  \phi \composition \psi$.
\end{theoremproof}

\begin{remark} \label{antichains:ergodicinvariantactions:remark}
  By Remark \ref{stratification:ergodicinvariantactions:remark}, after throwing out an $E$-invariant
  $\mu$-null \Borel set, the conclusion of Theorem \ref{antichains:ergodicinvariantactions} follows 
  from the weaker assumption that there is a \Borel subequivalence relation $F$ of $E$, for which
  $\cost{\mu}{F} > 1$, with respect to which $\mu$ is ergodic.
\end{remark}

\begin{remark}
  In particular, this gives a simple new proof of the existence of such actions for the
  equivalence relation generated by the usual action of $\SL{2}{\Z}$ on $\T[2]$, a result
  which originally appeared in \cite{Hjorth}.
\end{remark}

We next establish an analogous result in the absence of ergodicity.

\begin{theorem} \label{antichains:invariantactions}
  Suppose that $X$ is a \Polish space, $E$ is a projectively separable countable
  \Borel equivalence relation on $X$, $\mu$ is an $E$-invariant \Borel probability measure
  on $X$, and there is a free \Borel action of a countable non-abelian free group on $X$
  generating a subequivalence relation of $E$. Then for every countable non-abelian free
  group $\Gamma$, there is a \Borel sequence $\sequence{\cdot_r}[r \in \R]$ of free
  actions of $\Gamma$ on $X$, generating subequivalence relations $E_r$ of $E$, with the
  further property that $\sequence{E_r}[r \in \R]$ is an increasing sequence of relations which
  are pairwise incomparable under $\mu$-somewhere reducibility.
\end{theorem}

\begin{theoremproof}
  The proof is essentially the same as that of Theorem \ref
  {antichains:ergodicinvariantactions}. One minor difference is that
  we use Proposition \ref{stratification:invariantactions}
  in place of Proposition \ref{stratification:ergodicinvariantactions}
  to obtain a \Borel $\mu$-stratification $\sequence{\cdot_r}[r \in \R]$
  of $E$ by free actions of $\Gamma$ on $X$. Another difference is
  that we use $R$ to denote the relation on $\R$ in which two
  real numbers $r$ and $s$ are related if $E_r$ is merely $\mu$-somewhere
  reducible to $E_s$. As $\mu$-somewhere reducibility is weaker than
  $\mu$-reducibility in the absence of ergodicity, we must be slightly
  more careful in establishing the analog of Lemma \ref
  {antichains:ergodicinvariantactions:sections}.

  \begin{lemma}
    Every horizontal section of $R$ is countable.
  \end{lemma}

  \begin{lemmaproof}
    Suppose, towards a contradiction, that there exists $t \in \R$ for which $\horizontalsection
    {R}{t}$ is uncountable. For each $s \in \horizontalsection{R}{t}$, fix a $\mu$-positive \Borel
    set $B_s \subseteq X$ on which there is a \Borel reduction $\phi_s \from B_s \to X$ of
    $E_s$ to $E_t$. Then there exists $\epsilon > 0$ with $\mu(B_s) \ge \epsilon$ for
    uncountably many $s \in \horizontalsection{R}{t}$. As each $\phi_s$ is a homomorphism
    from $\restriction{(\intersection[r \in \R][E_r])}{B_s}$ to $E$, the $\mu$-nowhere
    hyperfiniteness of $\intersection[r \in \R][E_r]$ and the projective separability of
    $E$ ensure the existence of distinct $r, s \in \horizontalsection{R}{t}$ for which $\mu
    (B_r), \mu(B_s) \ge \epsilon$ and $\uniformmetric{\mu}(\phi_r, \phi_s) < \epsilon$. Then
    $\set{x \in B_r \intersection B_s}[\phi_r(x) = \phi_s(x)]$ is a $\mu$-positive \Borel set on
    which $E_r$ and $E_s$ coincide, a contradiction.
  \end{lemmaproof}

  One can now proceed exactly as in the proof of Theorem \ref{antichains:ergodicinvariantactions}
  to obtain the desired sequence.
\end{theoremproof}

\begin{remark} \label{antichains:invariantactions:remark}
  By Remark \ref{stratification:invariantactions:remark}, after throwing out an $E$-invariant
  $\mu$-null \Borel set, the conclusion of Theorem \ref{antichains:invariantactions} follows
  from the weaker assumption that for all $\mu$-positive \Borel sets $B \subseteq X$, there
  is a \Borel subequivalence relation $F$ of $\restriction{E}{B}$ for which $\cost{\restriction
  {\mu}{B}}{F} > \mu(B)$.
\end{remark}

Finally, we establish an analogous result in the absence of invariance.

\begin{theorem} \label{antichains:actions}
  Suppose that $X$ is a \Polish space, $E$ is a projectively separable countable
  \Borel equivalence relation on $X$, $\mu$ is a \Borel probability measure
  on $X$, and there is a free \Borel action of a countable non-abelian free group on $X$
  generating a $\mu$-nowhere hyperfinite subequivalence relation of $E$. Then for every
  countable non-abelian free group $\Gamma$, there is a \Borel sequence $\sequence
  {\cdot_r}[r \in \R]$ of free actions of $\Gamma$ on $X$, generating subequivalence
  relations $E_r$ of $E$, with the further property that $\sequence{E_r}[r \in \R]$ is an
  increasing sequence of relations which are pairwise incomparable under $\mu$-somewhere
  reducibility.
\end{theorem}

\begin{theoremproof}
  This follows from the proof of Theorem \ref{antichains:invariantactions}, using
  Proposition \ref{stratification:actions} in place of Proposition \ref
  {stratification:ergodicinvariantactions}.
\end{theoremproof}

\begin{remark} \label{antichains:actions:remark}
  By Remark \ref{stratification:actions:remark}, after throwing out an $E$-invariant
  $\mu$-null \Borel set, the conclusion of Proposition \ref{stratification:actions} follows from
  the weaker assumption that for all $\mu$-positive \Borel subsets $B \subseteq X$ such
  that $\restriction{\mu}{B}$ is equivalent to an $E$-invariant \Borel probability measure $\nu$,
  there is a \Borel subequivalence relation $F$ of $\restriction{E}{B}$ for which
  $\cost{\restriction{\nu}{B}}{F} > \nu(B)$.
\end{remark}

\begin{acknowledgments}
  We would like to thank the referee for his many useful suggestions.
\end{acknowledgments}

\bibliographystyle{amsalpha}
\bibliography{bibliography}

\end{document}